\documentclass[a4paper]{article}

\usepackage{amsmath,amssymb}

\newtheorem{thm}{Theorem}
\newtheorem{prop}[thm]{Proposition}

\newtheorem{cor}[thm]{Corollary}

\date{}

\begin{document}

\title{
Note on the Type II codes of length $24$
}

\author{S. Nagaoka and M. Oura}

\maketitle

\begin{abstract}
We express the weight enumerators of self-dual and doubly even (Type II for short) codes of length $24$
with a specified basis.
As a consequence, we present some congruence relations among the weight enumerators.
\end{abstract}

Keywords: code, weight enumerator.

2020MSC: Primary 95B05, Secondary 05E15.

\section{Introduction}
The weight enumerators of codes play important roles in coding theory.
Here we mainly focus on the properties of the weight enumerators
rather than on its applications.

On the other hand, the correspondences between codes and lattices are known.
In the paper \cite{NagaokaTakemori},
they studied the theta series of even unimodular lattices of length $24$.
Among other results,
they show that for an even unimodular lattice $\mathcal{L}$ with Coxeter number $h$,
we have
\begin{gather*}
\vartheta_{\mathcal{L}}^{(3)}=
(E_4^{(3)})^3+24(h-30)Y_{12}^{(3)}+48(h-30)^2X_{12}^{(3)}\\
+24(h-30)(2h^2+48h+1571)F_{12}
\end{gather*}
where $E_4^{(3)},Y_{12}^{(3)},X_{12}^{(3)},F_{12}$ are suitable Siegel modular
forms in genus $3$ with integral Fourier coefficients
and in particular $F_{12}$ vanishes under the action of Siegel's Phi operator.
The purpose of this note is to give the similar results in coding theory.

\section{Preliminaries}
\noindent

Let $\mathbf{F}_2=\{0,1\}$ be the field of two elements.
We sometimes regard as $\mathbf{F}_2\subset \mathbf{Z}$.
As is usual, we have the vector space $\mathbf{F}_2^n$
equipped with the inner product
\[
u\cdot v=u_1v_1+\dots + u_n v_n
\]
for $u=(u_1,\dots,u_n), \ y=(v_1,\dots,v_n)\in \mathbf{F}_2^n$.
The weight $wt(u)$ of $u\in \mathbf{F}_2^n$ is the number of non-zero coordinates of $u$.
A subspace $C$ of $\mathbf{F}_2^n$ is called a (linear) code of length $n$.
The dual code of $C$ is defined by 
\[
C^{\perp}=\{u\in \mathbf{F}_2^n:\ u\cdot v=0, \ \forall v \in C\}.
\]
If $wt(u)$ is a multiple of $4$ for all $u\in C$,
then $C$ is said to be doubly even.

We give some codes with generator matrices in which rows generate each code.
We put
\[
d_n:\
\begin{pmatrix}
11111 1& \dots & 0000\\
001111& \dots & 0000\\
& \ddots & \\
000000 &\dots  & 1111
\end{pmatrix}
\]
for $n\equiv 0\pmod{2}$ and 
\[
e_7:\ 
\begin{pmatrix}
0111100\\
0110011\\
1101010
\end{pmatrix},
\]

\[
e_8:\
\begin{pmatrix}
11110000\\
00111100\\
00001111\\
10101010
\end{pmatrix}.
\]
We denote by $g_{24}$ the binary Golay code of length $24$.

In this note, we deal with the self-dual and doubly even codes, Type II codes for short.
It is known that a Type II code of length $n$ exists 
if and only if $n\equiv 0\pmod{8}$.
Classification of Type II codes is completed up to $n=40$, 
see \cite{P1972,ps1975,CPS1992,BHM}.
Type II codes of length  $24$ are presented at Table \ref{tab:TypeIIcodes24}.
The $7$th code in that table is the binary Golay code $g_{24}$.

\bigskip

\begingroup
\renewcommand{\arraystretch}{1.5}
\begin{table}
\caption{Classification of Type II codes of length $24$.}\label{tab:TypeIIcodes24}
\[
\begin{tabular}{|c|c|c|c|c|c|c|c|c|c|}
\hline
$i$ & 1 & 2 & 3 &  4 &  5 & 6 & 7 & 8 &9 \\ \hline
Components  & $d_{12}^2$ & $d_{10}e_7^2$ & $d_8^3$ & $d_6^4$ & $d_{24}$ 
 & $d_4^6$ &  $g_{24}$   & $d_{16}e_8$ & $e_8^3$\\ \hline
 $h_i$ & $\frac{5}{4}$ & $1$ & $\frac{3}{4}$  & $\frac{1}{2}$ & $\frac{11}{4}$ & 
 $\frac{1}{4}$  & $0$ & $\frac{7}{4}$ & $\frac{7}{4}$ \\ \hline
\end{tabular}
\]
\end{table}
\endgroup

\bigskip

\noindent
The number $h_i$ is obtained as the number of elements of weight $4$ of each component divided by the dimension.
This number can be read off from the weight enumerator given below, that is,
a coefficient of $x^{n-4}y^4$ divided by $n$.
\noindent

For a code $C$, the weight enumerator of $C$ in genus $g$  is defined by
\[
W_C^{(g)}(x_a : \ a \in \mathbf{F}_2^g)=\sum_{u_1,\dots,u_g\in C}
\prod_{a\in \mathbf{F}_2^g}x_a^{n_a(u_1,\dots,u_g)}
\]
where $n_a(u_1,\dots,u_g)=\sharp \{i: a=(u_{1i},u_{2i},\dots,u_{gi})\}$.
In genus $1$, we may use $x,y$ instead of $x_0,x_1$.
If $C$ is of length $n$, we have the usual weight enumerator
\[
W_C^{(1)}(x,y)=\sum_{u\in C}x^{n-wt(u)}y^{wt(u)}.
\]
In the following, we may write $W_C^{(g)}$ for simplicity.
For codes $C$ and $C'$, we have $W_{C\oplus C'}^{(g)}=
W_{C}^{(g)}W_{C'}^{(g)}$.

For a column vector $a\in \mathbf{F}_2^g$, we define a map 
\begin{align*}
\Phi : \mathbf{C}[x_a\in \mathbf{F}_2^g]& \rightarrow \mathbf{C}[x_{a'}:a' \in \mathbf{F}_2^{g-1}]\\
x_{a}& \mapsto \begin{cases} x_{a'} & \text{ if }a=\begin{pmatrix} a' \\ 0\end{pmatrix} ,\\
0 & \text{ if }a=\begin{pmatrix} a' \\ 1 \end{pmatrix}.\end{cases}
\end{align*}

\noindent
It holds
$\Phi(W_C^{(g)})=W_C^{(g-1)}$.

It is known that
the ring generated over $\mathbf{C}$ 
by the weight enumerators of Type II codes in genus $g$
coincides with the invariant ring of some finite group, 
see \cite{gleason, duke, Runge}.
In particular, a basis of the vector pace generated over $\mathbf{C}$ 
by the weight enumerators of Type II code of length $24$ in $g=1,2$ is 
\[
W_{C_9}^{(1)},W_{C_7}^{(1)}
\]
and
\[
W_{C_9}^{(2)},W_{C_7}^{(2)},W_{C_5}^{(2)}.
\]

For completeness, we add here the weight enumerators in genus $1$.
\begin{align*}
W_{d_4}^{(1)}&=x^4+y^4,\\
W_{d_6}^{(1)}&=x^6+3x^2y^4,\\
W_{e_7}^{(1)}&=x^7+7x^3y^4,\\
W_{d_8}^{(1)}&=x^8+6x^4y^4+y^8,\\
W_{e_8}^{(1)}&=x^8+14x^4y^4+y^8,\\
W_{d_{10}}^{(1)}&=x^{10}+10x^6y^4+5x^2y^8,\\
W_{d_{12}}^{(1)}&=x^{12}+15x^8y^4+15x^4y^8+y^{12},\\
W_{d_{16}}^{(1)}&=x^{16}+28x^{12}y^4+70x^8y^8+28x^4y^{12}+y^{16},\\
W_{d_{24}}^{(1)}&=x^{24}+66x^{20}y^4+495x^{16}y^8+924x^{12}y^{12}+495x^8y^{16}+66x^4y^{20}+y^{24}.
\end{align*}
The coefficient of $x^{n-4}y^4$ is the number of elements of weight $4$ and 
each component of a Type II code of length $24$ has the same number, that is, $nh_i$.
For the Type II codes of length $24$, we mention
\begin{align*}
W_{C_5}^{(1)}&=x^{24}+66x^{20}y^4+495x^{16}y^8+2972x^{12}y^{12}+495x^8y^{16}+66x^4y^{20}+y^{24},\\
W_{C_7}^{(1)}&=x^{24}+759x^{16}y^8+2576x^{12}y^{12}+759x^8y^{16}+y^{24}
\end{align*}
and
\[
W_{C_5}^{(1)}=\frac{11}{7}W_{C_9}^{(1)}-\frac{4}{7}W_{C_7}^{(1)}.
\]

\section{Results}

We start with the case $g=1$.
Let
$\Delta =x^4y^4(x^4-y^4)^4$.
The vector space of the weight numerators of Type II codes in genus $1$ is
spanned by $W_{C_9}^{(1)}$ and $\Delta$. By direct calculation, we get the following theorem.
\begin{thm}\label{thm:weg1}
\begin{enumerate}
\item[\rm{(1)}]
For $i=1,2,\dots,9$, we have 
\[
W_{C_i}^{(1)}=W_{C_9}^{(1)}+6(4h_i-7)\Delta.
\]

\item[\rm{(2)}]
Let $i$ and $j$ be distinct integers in $1,2,\dots,8$.
If $4h_i\equiv 4h_j\pmod{m}$ for an integer $m$, then
\[
W_{C_i}^{(1)}\equiv W_{C_j}^{(1)}\pmod{6m}.
\]

\item[\rm{(3)}]
Let $C_{\alpha},C_{\beta}$ be Type II codes of length $24$ with $h_{\alpha}<h_{\beta}$.
Then we have
\[
W_{C_i}^{(1)}=\frac{h_i-h_{\beta}}{h_{\alpha}-h_{\beta}}W_{C_{\alpha}}^{(1)}+
\frac{h_i-h_{\alpha}}{h_{\beta}-h_{\alpha}}W_{C_{\beta}}^{(1)}
\]
for $i=1,2,\dots,9$.
\end{enumerate}
\end{thm}

\textit{Proof}.
We need only to prove (3). 
Set
\[
W_{C_i}^{(1)}=aW_{C_{\alpha}}^{(1)}+bW_{C_{\beta}}^{(1)}.
\]
Applying (1), we get a system of equations
\[
\begin{cases}
a+b=1,\\
6(4h_{\alpha}-7)a+6(4h_{\beta}-7)b=6(4h_i-7).
\end{cases}
\]
Since the determinant of the matrix 
$\begin{pmatrix} 1 & 1 \\ 6(4h_{\alpha}-7) & 6(4h_{\beta}-7)\end{pmatrix}$
is 
$-24(h_{\alpha}-h_{\beta})\not= 0$,
we get $a$ and $b$.
This completes the proof of Theorem \ref{thm:weg1}.

\bigskip

We add a few words on (2) of Theorem \ref{thm:weg1}.
For $i<j$, we denote by $m$ the number presented at the $(h_i,h_j)$-entry in Table \ref{tab:m}.
Then  (2) of Theorem \ref{thm:weg1} says that $(W_{C_i}^{(1)}-W_{C_j}^{(1)})/6m$ is in $\mathbf{Z}[x,y]$.
One can say more.
By direct calculation, we observe that 
the resulting $(W_{C_i}^{(1)}-W_{C_j}^{(1)})/6m$ contains a monomial with
coefficient $ 1$ or $-1$.

\begin{table}
\caption{Possible $m$ in (2) of Theorem\ref{thm:weg1}.}\label{tab:m}
\[
\begin{tabular}{|c|c|c|c|c|c|c|c|c|}
\hline
&  $h_2$ & $h_3$ & $h_4$ & $h_5$ & $h_6$& $h_7$ & $h_8$ \\
\hline
$h_1$  & 1 & 2 &3 &6 &4 & 5&2 \\
\hline
$h_2$  & &1 &2 &7 & 3& 4& 3 \\
\hline
 $h_3$ & & &1 &8 &2 &3 & 4 \\
\hline
 $h_4$ & & & &9 & 1& 2& 5 \\
\hline
 $h_5$ & & & & &10 &11 & 4\\
\hline
$h_6$  & & & & & &1 &6  \\
\hline
$h_7$  & & & & & & &7  \\
\hline
\end{tabular}
\]
\end{table}

\bigskip

We consider the case $g=2$.
The vector space of the weight enumerators of 
Type II codes of length $24$ in genus $2$ is 
spanned by $W_{C_9}^{(2)},W_{C_7}^{(2)},W_{C_5}^{(2)}$.
Let
\[
F=aW_{C_9}^{(2)}+bW_{C_7}^{(2)}+cW_{C_5}^{(2)}.
\]
Now we consider the action of $\Phi$ on $F$.
In order to make our discussion smooth, we set
\begin{align*}
X&=\frac{1}{42}\left(W_{C_9}^{(2)}-W_{C_7}^{(2)}\right),\\
Y&=-\frac{11}{7}W_{C_9}^{(2)}+\frac{4}{7}W_{C_7}^{(2)}+W_{C_5}^{(2)}.
\end{align*}
Since $\Phi(W_C^{(g)})=W_C^{(g-1)}$, we have
\begin{align*}
\Phi(F)&=aW_{C_9}^{(1)}+bW_{C_7}^{(1)}+cW_{C_5}^{(1)}\\
&=aW_{C_9}^{(1)}+bW_{C_7}^{(1)}+c\left( \frac{11}{7}W_{C_9}^{(1)}-
\frac{4}{7}W_{C_7}^{(1)}\right)\\
&= \left(a+\frac{11}{7}c\right)W_{C_9}^{(1)}+\left(b-\frac{4}{7}c\right)W_{C_7}^{(1)}.
\end{align*}
Then 
\begin{align*}
\Phi(F)=0&\Leftrightarrow a+\frac{11}{7}=b-\frac{4}{7}=0\\
&\Leftrightarrow a=-\frac{11}{7}c, b=\frac{4}{7}c\\
&\Leftrightarrow F=cY.
\end{align*}
Also
\begin{align*}
\Phi(F)=x^4y^4(x^4-y^4)^4
&\Leftrightarrow a+\frac{11}{7}=\frac{1}{42}, b-\frac{4}{7}=-\frac{1}{42}\\
&\Leftrightarrow a=\frac{1}{42}-\frac{11}{7}c, b=-\frac{1}{42}+\frac{4}{7}c\\
&\Leftrightarrow F=X+cY.
\end{align*}
We have thus obtained the following proposition.
\begin{prop}\label{prop:phi}
\begin{enumerate}
\item[\rm{(1)}]
$\Phi(F)=0$ if and only if $ F=cY$ for some constant $c$.

\item[\rm{(2)}]
$\Phi(F)=\Delta$ if and only if $F=X+cY$ for some constant $c$.
\end{enumerate}
\end{prop}

We introduce the following polynomials of
$\mathbf{Z}[x_a:a \in \mathbf{F}_2^2]$ from \cite{Oura2008}:

\begin{align*}
X_{24}&=X-\frac{1}{44}Y, \\
 Y_{24}&=\frac{1}{2^43\cdot 11}Y.
 \end{align*}
 
\noindent
Polynomials $W_{C_9}, X_{24}, Y_{24}$ form a basis of the vector space of
the weight enumerators of Type II codes of length $24$.
By Proposition \ref{prop:phi}, we have
\begin{align*}
\Phi(X_{24})&=\Phi(X)-\frac{1}{44}\Phi(Y)\\
&=\Delta
\end{align*}
and
\begin{align*}
\Phi(Y_{24})&=\frac{1}{2^43\cdot 11}\Phi(Y)\\
&=0.
\end{align*}
The coefficient $c(h_i)$ of $Y_{24}$ in 
\[
W_{C_i}=W_{C_9}^{(2)}+6(4h_i-7)X_{24}+c(h_i)Y_{24}
\]
can be obtained by direct calculation.
Therefore we have the following theorem 
corresponding to that of theta series in \cite{NagaokaTakemori}
mentioned in Introduction.

\begin{thm}\label{thm:weg2}
\begin{enumerate}
\item[\rm{(1)}]
$\Phi(X_{24})=\Delta$ and $\Phi(Y_{24})=0$.

\item[\rm{(2)}]
For $i=1,2,\dots,9$, we have
\begin{align*}
W_{C_i}^{(2)}&=W_{C_9}^{(2)}+6(4h_i-7)X_{24}+24(2h_i+3)(4h_i-7)Y_{24}.
\end{align*}
\end{enumerate}
\end{thm}


\begin{cor}
Let $i$ and $j$ be distinct integers in $1,2,\dots,8$.
If $4h_i\equiv 4h_j\pmod{m}$ for an integer $m$,
then 
\[
W_{C_i}^{(2)}\equiv W_{C_j}^{(2)} \pmod{6m}.
\]
\end{cor}

In order to state the next corollary, we introduce the Lagrange polynomial $\ell_{\epsilon}(x)$.
For $\epsilon \in \{\alpha,\beta,\gamma\}$, we set
\[
\ell_{\epsilon}(x)=\prod_{\substack{\mu \in \{\alpha,\beta,\gamma\}\\ \mu \not= \epsilon}}
\frac{x-x_{\mu}}{x_{\epsilon}-x_{\mu}}.
\]

\begin{cor}\label{cor:lagrange}
Let $C_{\alpha},C_{\beta},C_{\gamma}$ be Type II codes of length $24$
with $h_{\alpha}<h_{\beta}<h_{\gamma}$.
Then for $i=1,2,\dots,9$, the weight enumerator $W_{C_i}^{(2)}$ has 
the following expression:
\[
W_{C_i}^{(2)}=\ell_{\alpha}(h_i)W_{C_{\alpha}}^{(2)}+
\ell_{\beta}(h_i)W_{C_{\beta}}^{(2)}+
\ell_{\gamma}(h_i)W_{C_{\gamma}}^{(2)}.
\]
\end{cor}

\textit{Proof}. 
We set
\[
W_{C_i}^{(2)}=aW_{C_{\alpha}}^{(2)}+
bW_{C_{\beta}}^{(2)}+
cW_{C_{\gamma}}^{(2)}.
\]
We apply the expression in Theorem \ref{thm:weg2} to each $W_{C_{\alpha}}^{(2)},
W_{C_{\beta}}^{(2)},
W_{C_{\gamma}}^{(2)}$,
we get a system of equations
\[
A
\begin{pmatrix} a \\ b \\ c \end{pmatrix}
=\begin{pmatrix}
1 \\ c_0(h_i)\\ c_1(h_i)
\end{pmatrix}
\]
where 
\[
A=\begin{pmatrix} 1 & 1 & 1 \\
c_0(h_{\alpha}) & c_0(h_{\beta}) & c_0(h_{\gamma})\\
c_1(h_{\alpha}) & c_1(h_{\beta}) & c_1(h_{\gamma})
\end{pmatrix}
\]
and
\begin{align*}
c_0(h)&=6(4h-7),\\
c_1(h)&=24(2h+3)(4h-7).
\end{align*}
Since 
$\det A=-
4608(h_{\alpha}-h_{\beta})(h_{\alpha}-h_{\gamma})
(h_{\beta}-h_{\gamma})\not= 0$, 
we can solve the system of equations and get the result.
This completes the proof of Corollary \ref{cor:lagrange}.

\bigskip

There are two classes of Type II codes of length $16$.
A remark is that their weight enumerators are distinct in $g=3$.
This remark is important in number theory and we only mention a reference \cite{RungeSiegel}.
Inequality of the mentioned weight enumerators in $g=3$ leads to $W_{C_8}^{(3)}\not= W_{C_9}^{(3)}$.
Combining this with $h_8=h_9$, 
we see that Theorem \ref{thm:weg2} (2) can not be extended 
to the case $g=3$.

\bigskip

\textbf{Acknowledgment}.
Calculations are done with the help of Maple and Magma \cite{Bosmaetal}.
We thank Prof. Munemasa for comments to this manuscript.
The first named author is supported by JSPS KAKENHI Grant Number 20K03547.

\noindent
Shoyu Nagaoka\\
Yamato University, Japan\\
Emeritus Professor, Kindai University, Japan\\
shoyu1122.sn@gmail.com\\

\noindent
Manabu Oura\\
Kanazawa University, Japan\\
oura@se.kanazawa-u.ac.jp

\end{document}